\documentclass[12pt]{amsart}

\usepackage{epsfig}

\usepackage{amssymb}
\usepackage{amsmath}

\def\ds{\displaystyle}
\def\bs{\boldsymbol}

\newcommand{\N}{{\mathbb{N}}}
\newcommand{\R}{{\mathbb{R}}}

\newcommand{\SF}{{\mathbb{S}}}

\def\tl{\widetilde}
\def\ol{\overline}
\def\sm{\setminus}

\def\bu{{\mathbf{u}}}
\def\bv{{\mathbf{v}}}

\def\A{{\mathbf{A}}}
\def\B{{\mathbf{B}}}
\def\CM{\bs{\mathfrak{C}}}
\def\C{\mathcal{C}}

\def\ER{\mathbf{E}}
\def\BR{\mathbf{B}}
\def\I{{\mathbf{Id}}}
\def\K{\bs{\mathfrak{K}}}

\def\P{\mathcal{P}}
\def\LL{{\bs{\mathfrak{L}}}}
\def\NP{{\bs{\mathfrak{N}}}}

\makeatletter
\def\theequation{\thesection.\@arabic\c@equation}
\def\thethm{\@arabic\c@thm}

\makeatother

\newtheorem{theorem}{Theorem}[section]	
\newtheorem{lemma}[theorem]{Lemma}	       
\newtheorem{corollary}[theorem]{Corollary}	
\newtheorem{definition}[theorem]{Definition}	
\newtheorem{remark}[theorem]{Remark}	       

\title[Parabolic Systems with Discontinuous Data]%
   {A'priori Estimates and Precise Regularity for Parabolic Systems with
Discontinuous Data}

\author{D.K. Palagachev \and L.G. Softova}
\address{Dian K. Palagachev\\
   Politecnico di Bari,
   Dipartimento di Matematica,
   Via E. Orabona, 4,
   70~125~Bari,\ Italy}
   \email{dian@dm.uniba.it}
\address{Lubomira G. Softova\\
   Bulgarian Academy of Sciences,
   Institute of Mathematics,
   Sofia,\ Bulgaria}
   \email{luba@dm.uniba.it}


\begin{document}
\maketitle

\thispagestyle{empty}

\begin{abstract}
We deal with linear parabolic (in sense of Petrovskii) systems of order
$2b$ with discontinuous principal coefficients. A'priori estimates in Sobolev
and Sobolev--Morrey spaces are proved for the strong solutions by means of
potential analysis and boundedness of certain singular integral operators
with kernels of mixed homogeneity. As a byproduct, precise characterization
of the Morrey, $BMO$ and H\"older regularity is given for the solutions
and their derivatives up to order $2b-1.$
\end{abstract}

\section{Introduction}\label{sect1}
\setcounter{equation}{0}

The present  paper deals with linear systems of order $2b$ which are parabolic in the
sense of Petrovskii.  The
discontinuity of the principal coefficients $a_\alpha^{kj}(x,t)$ is expressed in terms
 of appurtenance of $a_\alpha^{kj}$'s to the
class of functions with {\it vanishing mean oscillation\/} which contains as a
proper subset the space of uniformly continuous functions.
We deal with strong solutions  belonging
to the Sobolev class $W^{2b,1}_{p,\mathrm{loc}}(Q),$ $p>1,$ or to
Sobolev--Morrey's space $W^{2b,1}_{p,\lambda,\mathrm{loc}}(Q)$ where $Q$ is a
cylindrical domain of $n$-dimensional base $\Omega$ and height $T.$ (Let us emphasize
at the very beginning that throughout the paper {\it ``local''\/} means
``local in space variables $x$ but global in time $t$''.) It is proved
that $a_\alpha^{kj}(x,t)\in VMO\cap L^\infty$ is a {\it sufficient\/} condition
ensuring local H\"older regularity of such  solutions and all their spatial derivatives
up to order $2b-1.$

Our approach makes use of the Calder\'on--Zygmund method of expressing highest order
 derivatives of the solution in terms of Gaussian-type potentials.
These turn out to be singular integrals with {\it kernels of mixed homogeneity\/}
of degree $-n-2b$ (strongly defined by the system itself) and their commutators
with the multiplication by the $VMO$ functions $a_\alpha^{kj}$ which have {\it small
integral oscillation\/} over small cylinders (cf. \cite{CFL} for the $L^p$-theory of
second-order elliptic equations,  \cite{CFF} and \cite{PS}
for the case of elliptic systems and higher-order elliptic equations, respectively,
and \cite{BC1,PS,Sf} for what concerns linear parabolic equations of second order).
Employing results on boundedness of these singular integral operators in $L^p$
(mainly due to Fabes--Rivi\`ere \cite{FR}) and in $L^{p,\lambda}$ (recently proven by
the authors in \cite{PS}), we derive a'priori bounds of Caccioppoli type
which yield estimates for the strong solutions in $W^{2b,1}_{p,\mathrm{loc}}(Q)$ (or
$W^{2b,1}_{p,\lambda,\mathrm{loc}}(Q)$) by means of the $L^{p}_{\mathrm{loc}}(Q)$
(respectively, $L^{p,\lambda}_{\mathrm{loc}}(Q)$) norm of the right-hand side plus a
weaker norm of the solution itself. By virtue of embedding properties of Sobolev
(Sobolev--Morrey) spaces into H\"older ones, these a'priori bounds lead to a complete
characterization of the Morrey, $BMO$ and H\"older regularity of the solution
and its spatial derivatives up to order $2b-1.$

The paper is organized as follows. In Section~\ref{sect2} we outline
the functional frameworks and state the main results. A geometric characterization
of solution's regularity is given in terms  of $n,$ $b,$ $p$ and $\lambda.$
Section~\ref{sect3} deals with the fundamental solution of parabolic systems with constant
coefficients. We study various properties of that fundamental solution which conduct
to a representation formula for solution's derivatives of order $2b$ through singular
integral operators with kernels of Calder\'on--Zygmund type.
 As far as H\"older's regularity follows from
known embedding between Sobolev and Besov spaces in the case of
$W^{2b,1}_{p,\mathrm{loc}}(Q)$ solutions, it is a rather delicate matter to deal with
solutions lying in $W^{2b,1}_{p,\lambda,\mathrm{loc}}(Q).$ For, we prove a
Poincar\'e-type inequality (Lemma~\ref{lem5}) which allows to employ recent
results (cf. \cite{K}) on Morrey, Campanato and H\"older classes on spaces of
homogeneous type. As outgrowth, we get a complete description of Morrey, $BMO$ and
H\"older regularity of the lower-order derivatives. Section~\ref{sect5}
contains various remarks regarding systems with lower order terms, non-zero initial
conditions and non-parabolic systems. We show that the parabolicity condition (see
\eqref{parabolicity} below) is not only sufficient but, in some context, also
necessary for the validity of the a'priori estimates obtained. Finally, a counter-example
is given to statements on H\"older's continuity that have appeared recently
in the literature.

A complete version containing complete proves of the presented here results is submitted.

\section{Main Results}\label{sect2}
\setcounter{equation}{0}

Let $\Omega\subset\R^n,$ $n\geq2,$ be a domain (open connected set) and
define $Q=\Omega\times (0,T)$ with $T>0.$
We consider the  linear system
\begin{equation}\label{system}
\LL(x,t,D_t,D^\alpha) \bu:= D_t \bu(x,t) - \sum_{|\alpha|=2b}\A_\alpha (x,t) D^\alpha \bu(x,t) =
 \mathbf{f}(x,t)
\end{equation}
for the unknown vector-valued function $\bu\colon Q\to\R^m$ given by the
transpose $\bu(x,t)=\big(u_1(x,t),\ldots,u_m(x,t)\big)^{\mathrm{T}},$
$\mathbf{f}=(f_1,\ldots,f_m)^{\mathrm{T}},$ and  $\A_\alpha(x,t)$
stands for the $m\times m$~matrix $\left\{ a_\alpha^{kj}(x,t) \right\}_{k,j=1}^m$
of the measurable coefficients $a_\alpha^{kj}\colon Q\to\R.$
Throughout the paper
$b\geq 1$ is a fixed integer, $\alpha=(\alpha_1,\ldots,\alpha_n)$ is a
multiindex of length $|\alpha|=\alpha_1+\cdots+\alpha_n,$
$D_t:=\partial/\partial t$  and
$D^\alpha\equiv D_x^\alpha:=D_1^{\alpha_1}\ldots
D_n^{\alpha_n}$ with $D_i:=\partial/\partial x_i.$
Further, $D^\alpha \bu=\big(D^\alpha u_1,\ldots,D^\alpha u_m\big)^{\mathrm{T}}$
and $D^s\bu$ substitutes {\it any\/} derivative $D^\alpha\bu$ with $|\alpha|=s\in\N.$
The boldface small roman letters $\bu,$ $\bv,$ $\mathbf{f},$ $\mathbf{g},\ldots$
denote $m$-dimensional vectors
whereas boldface capital letters $\mathbf{A},$ $\bs{\Gamma},\ldots$ stand for $m\times m$-matrices.
The notation $|\cdot|$ is used to indicate the Euclidean norm in $\R^N$ and $N$ will be
clear from the context.

We assume that the system    \eqref{system}  is uniformly parabolic in the
sense of Petrovskii (see \cite{Ed1,Ed2,F,Sl1}).
Namely, the $p$-roots of the $m$-degree polynomial
\begin{equation}\label{parabolic}
\text{det\,} \Big\{p\,\I_m -\sum_{|\alpha|=2b} \A_\alpha(x,t)(i\xi)^\alpha
\Big\}=0\qquad (i=\sqrt{-1})
\end{equation}
satisfy, for some $\delta>0$ and all $s=1,\ldots,m,$ the inequality
\begin{equation}\label{parabolicity}
\text{Re\,}p_s(x,t,\xi)\leq -\delta |\xi|^{2b}\qquad \text{for a.a.}\ (x,t)\in
Q,\  \forall
\xi\in \R^n.
\end{equation}
Here $\I_m$ is the identity $m\times m$~matrix and $\xi^\alpha
:=\xi_1^{\alpha_1}\xi_2^{\alpha_2}\cdots\xi_n^{\alpha_n}.$ Indeed, for
fixed $(x,t)\in Q$ and $\xi\in\R^n,$ $p_s(x,t,\xi)$ are nothing else than the
eigenvalues of the $m\times m$~matrix $(-1)^b
\sum_{|\alpha|=2b}\A_\alpha(x,t)\xi^\alpha,$ and the
parabolicity condition \eqref{parabolicity} means these have negative real part.

Our goal is to obtain interior H\"older regularity of the strong solutions to
\eqref{system} as a byproduct of a'priori estimates in Sobolev and Sobolev--Morrey
spaces. Let us recall the definitions of these functional classes.
\begin{definition}\label{dSob}
The parabolic Sobolev space $W^{2b,1}_p(Q),$ $p\in (1,+\infty),$  is the
collection of $L^p(Q)$ functions $u\colon Q\to\R$ all of which
distribution derivatives $D_t u$ and $D^\alpha_x u$ with $|\alpha|\leq 2b,$
belong to
$L^p(Q).$ The norm in $W^{2b,1}_p(Q)$ is
\begin{align*}
\|u\|_{W^{2b,1}_p(Q)}:=&\ \|D_t u\|_{p;Q} +\sum_{s=0}^{2b} \|D^s u\|_{p;Q}\\
=&\ \|D_t u\|_{p;Q} +\sum_{s=0}^{2b} \sum_{|\alpha|=s} \|D^\alpha_x
u\|_{p;Q},\qquad
 \|\cdot\|_{p;Q}:=\left(\int_Q|\cdot|^pdxdt \right)^{1/p}.
\end{align*}
For the sake of brevity, the cross-product of $m$ copies of $L^p(Q)$ is denoted by
the same symbol. Thus, if
$\mathbf{u}=(u_1,\ldots,u_m)^{\mathrm{T}}$ is a vector-valued function,
$\mathbf{u}\in L^{p}(Q)$ means that $u_k\in L^p(Q)$ for all $k=1,\ldots,m,$ and
$\|\mathbf{u}\|_{p;Q}:=\sum_{k=1}^m \|{u_k}\|_{p;Q}.$

When dealing with localized versions of $W^{2b,1}_{p}$ we always mean  local
in spatial variables $x$ and global in time, that is, $u\in
W^{2b,1}_{p,\mathrm{loc}}(Q)$ if $u\in W^{2b,1}_{p}(\Omega'\times(0,T))$ for any
$\Omega'\Subset\Omega.$
\end{definition}

Endow $\R^{n+1}=\R^n_x\times\R_t$ with the
{\it parabolic\/} metric $\varrho(x,t)=\max\{|x|,|t|^{1/2b}\}.$
We will employ the system of
{\it parabolic cylinders\/}
\begin{equation}\label{pc}
\C_r(x_0,t_0):= B_r(x_0)\times	(t_0-r^{2b},t_0),\quad B_r(x_0):= \big\{x\in \R^{n}\colon |x-x_0|<r,\big\}.
\end{equation}
 Obviously, the Lebesgue measure
$|\C_r|$ is comparable to $r^{n+2b}.$
\begin{definition}\label{dMor}
Let $p\in(1,+\infty)$ and $\lambda\in(0,n+2b).$
The function $u\in L^p(Q)$ belongs to the parabolic Morrey space
$L^{p,\lambda}(Q)$ if
\[
\|u\|_{p,\lambda;Q}:=\left(\sup_{r>0}\frac{1}{r^\lambda}\int_{\C_r\cap Q}
|u(x,t)|^p dxdt \right)^{1/p} <\infty
\]
where $\C_r$ ranges in the set of  parabolic cylinders in $\R^{n+1}.$
The Sobolev--Morrey space $W^{2b,1}_{p,\lambda}(Q)$
 consists of all functions $u\in
W^{2b,1}_p(Q)$ with generalized derivatives $D_tu$ and $D^\alpha_xu,$
$|\alpha|\leq 2b,$ belonging to $L^{p,\lambda}(Q).$ The norm in
$W^{2b,1}_{p,\lambda}(Q)$ is given by
$\|u\|_{W^{2b,1}_{p,\lambda}(Q)}:=\|D_tu\|_{p,\lambda;Q}+
\sum_{s=0}^{2b}\sum_{|\alpha|=s}\|D^\alpha_x u\|_{p,\lambda;Q}.$
\end{definition}
We refer the reader to \cite{Cm,Cm2,Da,K,M}
for various properties of Morrey and Sobolev--Morrey spaces.
\begin{definition}\label{dBMO}
For a locally integrable function $f\colon\ \R^{n+1}\to\R$ define
\[
\eta_f(R):=\sup_{r\leq R} \frac{1}{|\C_r|}\int_{\C_r}
 |f(y,\tau)-f_{\C_r}|dyd\tau\quad \text{for every}\ R>0,
\]
where $\C_r$ is any parabolic cylinder and
$f_{\C_r}$ is the average $\frac{1}{|\C_r|}\int_{\C_r} f(y,\tau) dyd\tau.$
Then:
\begin{itemize}
\itemsep=1pt
\item[$\bullet$]  $f\in BMO$ {\em (bounded mean oscillation\/}, see
John--Nirenberg~{\rm\cite{JN}}$)$ if $\|f\|_*:=\sup_R \eta_f(R)<+\infty.$
$\|f\|_*$ is a norm in $BMO$ modulo constant functions under which $BMO$
is a Banach space.
\item[$\bullet$] $f\in VMO$ {\em (vanishing mean oscillation\/}, see
Sarason~{\rm\cite{S}}$)$ if
$f\in BMO$ and $\lim_{R\downarrow 0}\eta_f(R)=0.$
The quantity $\eta_f(R)$ is referred to as  $VMO$-modulus of $f.$
\end{itemize}
The spaces $BMO(Q)$ and $VMO(Q),$ and $\|\cdot\|_{*;Q}$
 are defined in a similar manner
taking $\C_r\cap Q$ instead of $\C_r$ above.
\end{definition}

The main results of this paper are as follows.
\begin{theorem}\label{th1}
Suppose \eqref{parabolicity},	$1<q\leq p<+\infty,$
$a_\alpha^{kj}\in VMO(Q)\cap L^\infty(Q),$
$\mathbf{f}\in L^p_{\mathrm{loc}}(Q)$ and let $\bu\in
W^{2b,1}_{q,{\mathrm{loc}}}(Q)$  be a strong solution of
\eqref{system} such that $\bu(x,0)=\mathbf{0}.$
Then the operator $\LL$	improves integrability, that is,
$\bu\in W^{2b,1}_{p,{\mathrm{loc}}}(Q),$ and for any
$Q'=\Omega'\times(0,T),$ $Q''=\Omega''\times(0,T),$
$\Omega'\Subset\Omega''\Subset\Omega,$
there is a constant $C$ depending on
$n,$ $p,$ $m,$ $b,$ $\delta,$ $\|a^{kj}_\alpha\|_{\infty;Q},$
$\eta_{a^{kj}_\alpha}$ and
 $\mathrm{dist\,}(\Omega',\partial\Omega'')$ such that
 \begin{equation}\label{9}
\|\bu\|_{W^{2b,1}_p(Q')}\leq
C\left(\|\mathbf{f}\|_{p;Q''}+\|\bu\|_{p;Q''}\right).
\end{equation}
\end{theorem}
We will show in Section~\ref{sect5} that the parabolicity condition
\eqref{parabolicity} is also {\it necessary\/} in order to be valid the estimate
\eqref{9}.

Since $W^{2b,1}_p(Q')$ is contained into the Besov space
$B^{\sigma,\sigma/2b}_{\infty,\infty}(Q')$ with $\sigma=2b-\frac{n+2b}{p}>0$
and $B^{\sigma,\sigma/2b}_{\infty,\infty}(Q')$  coincides with the H\"older
space $C^{\sigma,\sigma/2b}(Q')$ for {\it non-integer\/} $\sigma$
(see \cite{I,ISl}, \cite[Theorems~2.5, 2.7]{Sl2}), we get
\begin{corollary}\label{cr1}
In addition to the hypotheses of Theorem~$\ref{th1},$ suppose
$p>\frac{n+2b}{2b}.$ Then $\bu\in L^\infty(Q')$ and there is a constant
$C$ such that
\[
\|\bu\|_{\infty;Q'}\leq
C\left( \|\mathbf{f}\|_{p;Q''}+\|\bu\|_{p;Q''}\right).
\]
Moreover, the $x$-derivatives of $\bu$ are H\"older continuous for large values of
$p.$ Precisely,
\begin{itemize}
\item[$\bullet$] if $p\in\left(\frac{n+2b}{2b-s},\frac{n+2b}{2b-s-1}\!\right)$
for a fixed $s\in\{0,1,\ldots,2b-2\}$ then
$D^s\bu \in C^{\sigma_s,\sigma_s/2b}(Q')$ with $\sigma_s=2b-s-\frac{n+2b}{p};$
\item[$\bullet$] if $p\in(n+2b,+\infty)$ then $D^{2b-1}\bu \in
C^{\sigma_{2b-1},
\sigma_{2b-1}/2b}(Q')$ with $\sigma_{2b-1}=1-\frac{n+2b}{p}$
\end{itemize}
and in all cases
\[
\sup_{\underset{(x,t),\ (x',t')\in Q'}{(x,t)\neq(x',t')}}
\frac{\left|D^{s}\bu(x,t)-D^{s}\bu(x',t')\right|}{\left(|x-x'|+|t-t'|^{1/2b}
\right)^{\sigma_{{s}}}} \leq
C\left( \|\mathbf{f}\|_{p;Q''}+\|\bu\|_{p;Q''}\right)
\]
for $s\in\{0,1,\ldots,2b-1\}.$
\end{corollary}

Our next result provides improving-of-integrability property of $\LL$ and
a'priori estimates in Sobolev--Morrey spaces for solutions of
\eqref{system} with Morrey right-hand side.
\begin{theorem}\label{th3}
Suppose \eqref{parabolicity},	$1<q\leq p<+\infty,$ $\lambda\in(0,n+2b),$
$a_\alpha^{kj}\in VMO(Q)\cap L^\infty(Q),$
$\mathbf{f}\in L^{p,\lambda}_{\mathrm{loc}}(Q)$ and let $\bu\in
W^{2b,1}_{q,{\mathrm{loc}}}(Q)$  be a strong
solution of \eqref{system} such that $\bu(x,0)=\mathbf{0}.$
Then $\bu\in W^{2b,1}_{p,\lambda,{\mathrm{loc}}}(Q)$ and
\begin{equation}\label{9a}
\|\bu\|_{W^{2b,1}_{p, \lambda}(Q')}\leq
C\left(\|\mathbf{f}\|_{p,\lambda;Q''}+
\|\bu\|_{p,\lambda;Q''} \right)
\end{equation}
with a constant $C$ depending on the quantities listed in
Theorem~$\ref{th1}$ and on $\lambda$ in addition.
\end{theorem}

An essential step in the prove of the a'priori estimate \eqref{9a} is  obtaining of local Caccioppoli-type estimate
for the upper derivatives of the solution
\begin{equation}\label{6}
\|D^{2b}\bu\|_{p, \lambda;\C_{r/2}}\leq
C\left(\|\LL \bu\|_{p,\lambda;\C_r}+ Cr^{-2b}
\|\bu\|_{p,\lambda;\C_r} \right)
\end{equation}
which holds for $\|D_t\bu\|_{p, \lambda;\C_{r/2}}$ as well by virtue of the parabolic structure of \eqref{system}.

As consequence of \eqref{9a} we obtain precise characterization
of Morrey, $BMO$ and H\"older regularity of the derivatives $D^s\bu$
with $s\in\{0,1,\ldots,$ $2b-1\}.$ Precisely,
\begin{corollary}\label{cr4}
Under the hypotheses of Theorem~$\ref{th3}$
fix an $s\in\{0,1,\ldots,$ $2b-1\}.$ Then there is a constant
$C$ such that
\begin{itemize}
\item[$\bullet$] if $p\in \left(1,\frac{n+2b-\lambda}{2b-s}\right)$ then
$D^s\bu \in L^{p,(2b-s)p+\lambda}(Q')$  and
\[
\|D^s\bu\|_{p,(2b-s)p+\lambda;Q'} \leq
C\left(
\|\mathbf{f}\|_{p,\lambda;Q''}+\|\bu\|_{p,\lambda;Q''}\right);
\]
\item[$\bullet$] if $p=\frac{n+2b-\lambda}{2b-s}$ then $D^{s}\bu \in BMO (Q')$
and
\[
\|D^s\bu\|_{*;Q'} \leq
C\left(
\|\mathbf{f}\|_{p,\lambda;Q''}+\|\bu\|_{p,\lambda;Q''}\right);
\]
\item[$\bullet$] if
$p\in\left(\frac{n+2b-\lambda}{2b-s},\frac{n+2b-\lambda}{2b-s-
1}\right)$\footnote{In the case
$s=2b-1$ this inclusion rewrites naturally as $p\in (n+2b-\lambda,+\infty).$}
 then
$D^s\bu \in C^{\sigma_s,\sigma_s/2b}(Q')$ with
$\sigma_s=2b-s-\frac{n+2b-\lambda}{p}$	and
\[
\sup_{\underset{(x,t),\ (x',t')\in Q'}{(x,t)\neq(x',t')}}
\frac{\left|D^s\bu(x,t)-D^s\bu(x',t')\right|}{\left(|x-x'|+|t-t'|^{1/2b}\right
)^{\sigma_s}} \leq
C\left(
\|\mathbf{f}\|_{p,\lambda;Q''}+\|\bu\|_{p,\lambda;Q''}\right).
\]
\end{itemize}
\end{corollary}

\begin{figure}[hbt]
\includegraphics[scale=.6]{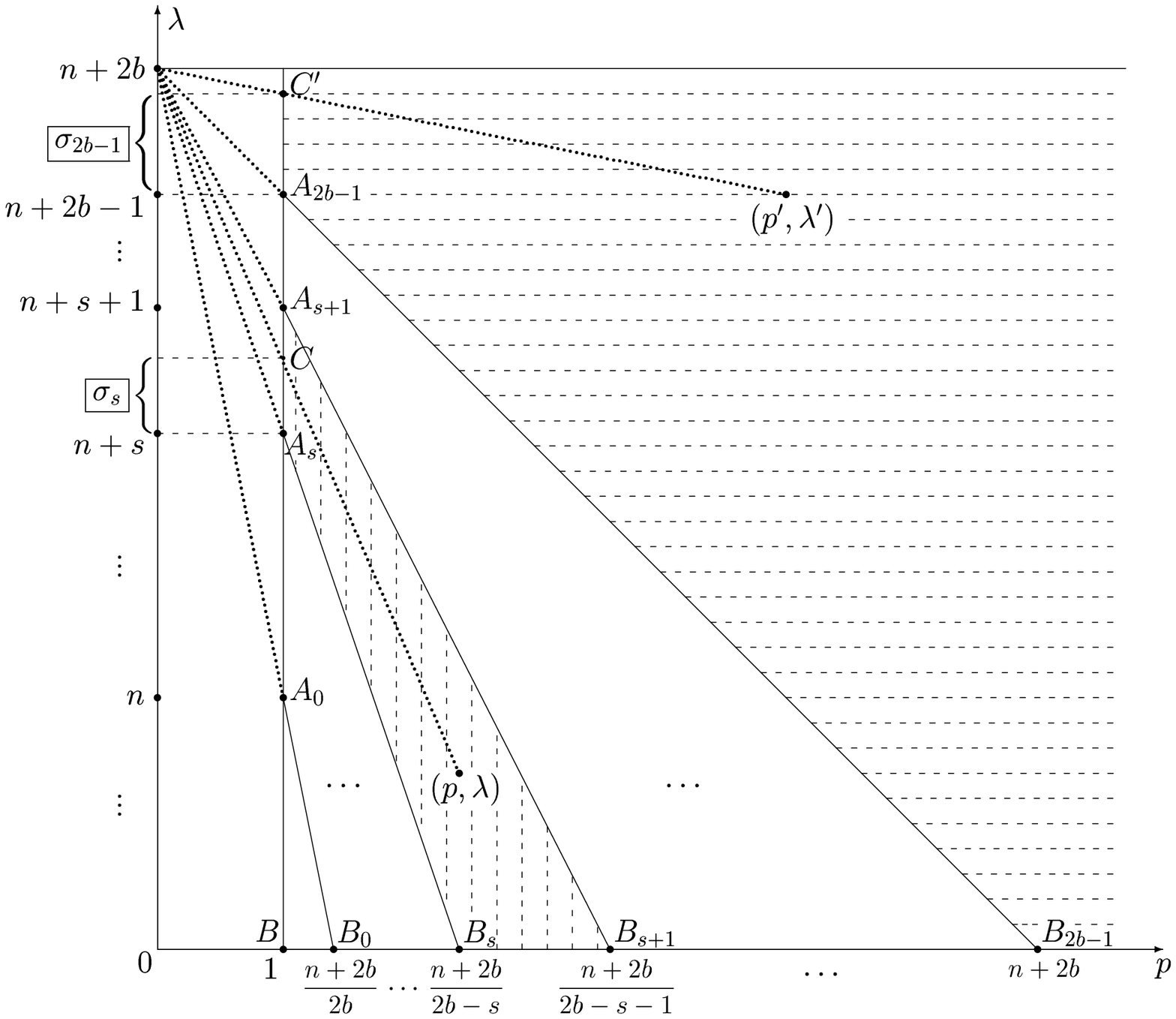}
\caption{The plane $O_{p,\lambda}$}
\label{picture}
\end{figure}
A simple geometric interpretation is proposed on Figure~$\ref{picture}$
of the results of Corollaries~$\ref{cr4}$ and $\ref{cr1}.$
A typical situation is considered for the couple $(p,\lambda)$
lying in the semistrip $\big\{(p,\lambda)\colon\ p>1,\ 0<\lambda<n+2b\big\}$
and $s\in\{0,1,\ldots,2b-1\}.$
The points $B_s$ on the $p$-axis are simply $\left(\frac{n+2b}{2b-s},0\right),$
$B=(1,0),$ and $A_s=(1,n+s)$ is the intersection of the line $\{p=1\}$ with the line
passing through $(0,n+2b)$ and $B_s.$

If $(p,\lambda)$ belongs to the open right triangle $BB_sA_s$ then $D^s\bu\in
L^{p,(2b-s)p+\lambda}(Q').$ In particular, $(p,\lambda)\in \triangle BB_0A_0$ yields
$\bu\in L^{p,2bp+\lambda}(Q')$ whereas $\bu\in BMO(Q')$ if
$(p,\lambda)$ lies on the open line segment $A_0B_0.$

Let $s\in\{0,\ldots,2b-2\}$ and
take $(p,\lambda)$ in the interior of the tetragon $R_s:=B_sB_{s+1}A_{s+1}A_s.$
Then the spatial derivatives $D^s\bu$ are H\"older continuous
with exponent $\sigma_s$ given by Corollary~$\ref{cr4},$ while $D^{s+1}\bu \in
L^{p,(2b-s-1)p+\lambda}(Q').$
Moreover, $\sigma_s$ is the length $|{CA_s}|$ of the segment ${CA_s}$ where
$C=C(p,\lambda)$ is the intersection of the vertical line $\{p=1\}$ with the
line connecting the points $(p,\lambda)$ and $(0,n+2b).$
If $(p,\lambda)\in A_sB_s$ (the open line segment)
we have  $D^s\bu \in BMO,$ whereas $(p,\lambda)\in A_{s+1}B_{s+1}$
implies $D^{s+1}\bu\in BMO.$

The situation is similar for the derivatives $D^{2b-1}\bu$ as well (i.e.,
$s=2b-1$) but $R_{2b-1}$ is now the shadowed {\it unbounded region\/} on the
picture. Thus, $(p',\lambda')\in R_{2b-1}$ gives that $D^{2b-1}\bu$ are H\"older
continuous with exponent $\sigma_{2b-1}=|{C'A_{2b-1}}|,$
$C'=C'(p',\lambda'),$ while $(p',\lambda')\in A_{2b-1}B_{2b-1}$ yields
$D^{2b-1}\bu\in BMO.$

To interpret the statement of Corollary~$\ref{cr1},$ we have simply to
consider values of $p>1$ lying on the $p$-axis. This way, $p>\frac{n+2b}{2b}$
implies $\bu \in L^\infty(Q')$ whereas $D^s\bu$ are H\"older continuous with
exponent $\sigma_s=2b-s-\frac{n+2b}{p}$ if $p$ lies on the open line segment
$B_sB_{s+1}$ (with the setting $B_{2b}:=+\infty$).

The proof of Corollary~\ref{cr4} relies on the next result
which is a parabolic version of the classical Poincar\'e inequality.
\begin{lemma}\label{lem5}
Let $\bu\in W^{2b,1}_p(\C_r),$ then for each $s\in\{0,1,\ldots,2b-1\}$ there is a constant
$C=C(p,m,n,s)$ such that
\begin{align*}
\int_{\C_r} |D^{s}\bu(x,t)-(D^{s}\bu)_{\C_r}|^p dxdt
\leq&\
C\big(r^{(2b-s)p}\left(\|D^{2b}\bu\|^p_{L^p(\C_r)}+\|D_t\bu\|^p_{L^p(\C_r)}
\right)\\
&\ +r^p\|D^{s+1}\bu\|^p_{L^p(\C_r)}\big).
\end{align*}
\end{lemma}

As  consequence of Lemma~\ref{lem5} and the absolute continuity of the Lebesgue integral
 we have the following	refinement of Corollary~$\ref{cr1}$
\begin{corollary}\label{cr5}
Under the hypotheses of Theorem~$\ref{th1},$
fix an $s\in\{0,1,\ldots,$ $2b-1\}.$ Then there is a constant
$C$ such that
\begin{itemize}
\item[$\bullet$] if $p\in \left(1,\frac{n+2b}{2b-s}\right)$ then
$D^s\bu \in L^{p,(2b-s)p}(Q')$  and
\[
\|D^s\bu\|_{p,(2b-s)p;Q'} \leq
C\left( \|\mathbf{f}\|_{p;Q''}+\|\bu\|_{p;Q''}\right);
\]
\item[$\bullet$] if $p=\frac{n+2b}{2b-s}$ then $D^{s}\bu \in BMO (Q'),$
\[
\|D^s\bu\|_{*;Q'} \leq
C\left( \|\mathbf{f}\|_{p;Q''}+\|\bu\|_{p;Q''}\right)
\]
and moreover, $D^{s}\bu \in VMO (Q').$
\end{itemize}
\end{corollary}
 Thus
$D^{s}\bu \in L^{p,(2b-s)p}(Q')$ if $p\in BB_s$ while $p=B_s$ gives
$D^s\bu\in VMO(Q').$

\section{Gaussian-type Potentials}\label{sect3}
\setcounter{equation}{0}

Consider now the metric $\rho(x,t)$ defined in \cite{FR}.
The mixed homogeneity of $\rho$ fits perfectly to the geometry of
$2b$-order parabolic operators and really, $\rho(x,t)$ was first introduced
by Fabes and Rivi\`ere(\cite{FR}) in the study of singular integral operators
with kernels of mixed homogeneity.
\begin{definition}\label{d1}  A function $k(x,t;y,\tau)\colon\
\R^{n+1}_{x,t}\times\big(\R^{n+1}_{y,\tau}\setminus \{0\}\big)\to \R$
is said to be  a {\em  variable Calder\'on-Zygmund kernel of parabolic type\/}
 if:
\begin{itemize}
\itemsep=3pt
\leftskip=8pt
\item[$ i)$] it is a  {\em constant Calder\'on-Zygmund kernel} for a.a.  fixed
$(x,t)\in \R^{n+1}.$   That is,
\begin{itemize}
\itemsep=3pt
\leftskip=16pt
\item[$i_a)$] $k(x,t;\cdot,\cdot)\in C^\infty(\R^{n+1}\setminus \{0\})$;
\item[$ i_b)$]	$k(x,t;\mu y_1,\ldots,
\mu y_n,\mu^{2b} \tau)=\mu^{-(n+2b) } k(x,t;y,\tau)\ $	for each $\mu>0;$
\item[$i_c)$] $\ds\int_{{\SF^{n}}} |k(x,t;y,\tau)|d\sigma_{(y,\tau)}<\infty$
and
$\ds\int_{{\SF^{n}}} k(x,t;y,\tau) d\sigma_{(y,\tau)} =0.$
\end{itemize}
\item[$ii)$]
$ \ds\sup_{(y,\tau)\in \SF^{n}}\left|D^\beta_{y,\,\tau}
 k (x,t;y,\tau)\right|\leq C(\beta)$ $\forall$ multiindex $\beta,$
independently of $(x,t).$
\end{itemize}
\end{definition}

Turning back to the system \eqref{system}, we fix $a_\alpha^{kj}$'s
at a point $(x_0,t_0),$ set
$$
\A_\alpha^0:=\A_\alpha(x_0,t_0)=\left\{ a_\alpha^{kj}(x_0,t_0)
\right\}_{k,j=1}^m,
$$
 and consider
the constant coefficients operator
\[
\LL_0(D_t,D^\alpha) :=\LL(x_0,t_0,D_t,D^\alpha) = \I_m D_t- \sum_{|\alpha|=2b}\A_\alpha^0 D^\alpha .
\]
It is known (see \cite{Sl1}) that the fundamental matrix
$\bs{\Gamma}_0(x,t)=\big\{
\Gamma_0^{kj}(x,t)\big\}_{k,j=1}^m$ of $\LL_0$ has entries
\[
\Gamma_0^{kj}(x,t)=L_{jk}(D_t,D^\alpha) \widetilde\Gamma_0(x,t),
\]
where $\{L_{jk}(D_t,D^\alpha)\}_{j,k=1}^m$ is the {\it cofactor matrix} of
$$
\LL_0(D_t,D^\alpha)=
\Big\{\delta^{jk}D_t-\sum_{|\alpha|=2b} a_\alpha^{jk}(x_0,t_0)D^\alpha
\Big\}_{j,k=1}^m,
$$
and $\widetilde\Gamma_0(x,t)$ is the fundamental solution of the {\it
parabolic equation}
\begin{equation}\label{det}
\text{det\,}\LL_0(D_t,D^\alpha) u =    \text{det\,}\Big\{\I_m D_t -
\sum_{|\alpha|=2b}\A_\alpha^0 D^\alpha\Big\}u=0.
\end{equation}
(Note that $L_{jk}$ is either a homogeneous differential operator of order $2b(m-1)$
or the operator of multiplication by $0.$)
Applying Fourier transform in $x$ and Laplace transform in $t,$ it is easy to
get
\begin{equation}\label{fund-equa}
\widetilde\Gamma_0(x,t) =\frac{1}{(2\pi)^n 2\pi i}\int_{\R^n} {\mathrm{e}}^{i(x\cdot\xi)}
d\xi \int_{{\mathfrak{C}}(\xi)}
\frac{{\mathrm{e}}^{pt}}{\text{det\,} \left\{p\,\I_m-\sum_{|\alpha|=2b} \A^0_\alpha
(i\xi)^\alpha \right\}} dp, \end{equation}
where ${\mathfrak{C}}(\xi)$ is a contour in the complex $p$-plane enclosing
all the roots of \eqref{parabolic}
and therefore, in view of \eqref{parabolicity}, could be taken to lie in the
left half-plane.

The fundamental matrix	$\bs{\Gamma}_0(x,t)$ possesses properties analogous to these
of the Gauss kernel (see \cite{Ed1,Ed2,F-ind,F,Sl1}). Precisely,
 \begin{itemize}
\leftskip=4pt
\item[$(\P_1)$] {\it Regularity:\/}  $\bs{\Gamma}_0\in C^\infty(\R^{n+1}\sm
\{0\}).$
\item[$(\P_2)$] {\it Mixed homogeneity:\/} for any $\mu>0$ and any
multiindex $\beta$ it holds
\[
\bs{\Gamma}_0(\mu x,\mu^{2b} t)=
\mu^{- n}\bs{\Gamma}_0(x,t),\qquad
D^\beta \bs{\Gamma}_0(\mu x,\mu^{2b} t)=
 \mu^{-n-|\beta|}D^\beta\bs{\Gamma}_0(x,t).
\]
\item[$(\P_3)$] {\it Vanishing property on the unit sphere $\SF^{n}:$\/}
\[
\int_{\SF^{n}} D^\alpha\bs{\Gamma}_0( x, t) d\sigma_{( x,t)}=\mathbf{0}\quad
\text{for any}\ \alpha,\ |\alpha|=2b.
\]
\item[$(\P_4)$] {\it Boundedness of the derivatives:\/}
\[
\sup_{(x, t)\in\SF^{n}}|D^\beta\bs{\Gamma}_0( x, t)|\leq C(n,\beta,
\max_{\alpha}|\A^0_\alpha|)
 \quad \forall\  \text{multiindex}\  \beta.
\]
\item[$(\P_5)$] {\it Integrability:\/}
\[
 D^\beta\bs{\Gamma}_0\in L^1_{\mathrm{loc}}(\R^{n+1}) \quad\text{for}\
|\beta|<2b,\qquad
D^\alpha\bs{\Gamma}_0\not\in L^1_{\mathrm{loc}}(\R^{n+1}) \quad\text{for}\
|\alpha|=2b.
\]
\end{itemize}

Take  an arbitrary function
$\bv\in C^\infty(\R^{n+1})$ which is compactly supported in $x$ and
$\bv(x,0)=\mathbf{0}.$ Let the point $(x_0,t_0),$ where the coefficient
matrix $\A_\alpha(x_0,t_0)$ were frozen, belong to $\text{supp\,}\bv.$
Making use of a standard approach and unfreezing the coefficients
(see \cite{BC1, LSU}) we obtain a representation formula for the derivatives of $\bv$ of order $2b$
\begin{align}\label{10}
D^\alpha \bv(x,t)=\  &\ p.v.
\int_{\R^{n+1}}D^\alpha\bs{\Gamma}(x,t;x-y,t-\tau)\LL{\bv}(y,\tau)dyd\tau\\
\nonumber
&\ +\sum_{|\alpha'|=2b} p.v. \int_{\R^{n+1}}D^{\alpha}\bs{\Gamma}(x,t;x-y,t-\tau)\\
\nonumber
&\qquad\qquad\qquad\qquad
\times	\big(\A_{\alpha'}(y,\tau)-
\A_{\alpha'}(x,t)\big) D^{\alpha'}_y
\bv(y,\tau)dyd\tau\\
\nonumber
&\ +   \int_{\SF^{n}}D^{\beta^s}\bs{\Gamma}(x,t;y,\tau)\nu_s
d\sigma_{(y,\tau)}\,\LL{\bv}(x,t)\\
\nonumber
 =:&\ \K_\alpha (\LL\mathbf{v})+\sum_{|\alpha'|=2b} \CM_\alpha[\A_{\alpha'},
D^{\alpha'} \bv]+
 \mathbf{F}(x,t) \LL{\bv}(x,t),\qquad \forall \alpha\colon |\alpha|=2b
 \end{align}
where the derivatives $D^\alpha\bs{\Gamma}(\cdot,\cdot;\cdot,\cdot)$
are taken with respect to the third variable.

Denote $\mathbf{k}(x,t;y,\tau):=D^\alpha_y\bs{\Gamma}(x,t;y,\tau)$ with
$|\alpha|=2b.$
Each entry of the $m\times m$~matrix $\mathbf{k}$ is a Calder\'on-Zygmund kernel
in the sense of the Definition~$\ref{d1}.$ In fact, $i_a)$ and $i_b)$ are
just properties $(\P_1)$ and $(\P_2)$ of the fundamental solution, while
$i_c)$ and $ii)$ follow from $(\P_3)$ and $(\P_4).$
Finally, $(\P_5)$ shows that $\K_\alpha$ and $\CM_\alpha$ are really  singular
integral operators. For what concerns their boundedness
in Lebesgue and Morrey spaces, we have
\begin{lemma}\label{pr1}
Let $|\alpha|=|\alpha'|=2b$ and $\A_\alpha\in L^\infty(Q).$
For each $p\in(1,+\infty)$ there exists a constant
$C=C(n,m,b,\delta,\|\A_\alpha\|_{\infty;Q},p)$ such that
\begin{align}\label{eqKf}
\|\K_\alpha \mathbf{f}\|_{p;Q}&\leq C \|\mathbf{f}\|_{p;Q},\\
\label{eqCaf}
\|\CM_\alpha [\A_{\alpha'}, \mathbf{f} ]\|_{p;Q}&\leq C \|\A_{\alpha'}\|_{*;Q}
\|\mathbf{f}\|_{p;Q}
\end{align}
for any $\mathbf{f}\in L^p(Q).$

For each $p\in(1,+\infty)$ and each $\lambda\in(0,n+2b)$ there is a constant
$C$
depending on $n,$ $m,$ $b,$ $\delta,$ $\|\A_\alpha\|_{\infty;Q},$ $p$ and
$\lambda$ such that
\begin{align}\label{eqKfMor}
\|\K_\alpha \mathbf{f}\|_{p,\lambda;Q}&\leq C
\|\mathbf{f}\|_{p,\lambda;Q},\\
\label{eqCafMor}
\|\CM_\alpha [\A_{\alpha'}, \mathbf{f} ]\|_{p,\lambda;Q}&\leq C \|\A_{\alpha'}\|_{*;Q}
\|\mathbf{f}\|_{p,\lambda;Q}
\end{align}
for any $\mathbf{f}\in L^{p,\lambda}(Q).$

Moreover, let $\A_\alpha\in VMO(Q)\cap L^\infty(Q)$ with $VMO$-modulus $\eta_{\A_\alpha}.$
Then for each $\varepsilon>0$ there exists $r_0=r_0(\varepsilon,\eta_{\A_\alpha})$
such
that if $r<r_0$ we have
\begin{align}
\label{3.10}
\|\CM_\alpha [\A_{\alpha'}, \mathbf{f} ]\|_{p;\C_r}&\leq C \varepsilon
\|\mathbf{f}\|_{p;\C_r}\qquad \forall \ \mathbf{f}\in L^p(\C_r),\\
\label{3.11}
\|\CM_\alpha [\A_{\alpha'}, \mathbf{f} ]\|_{p,\lambda;\C_r}&\leq C \varepsilon
\|\mathbf{f}\|_{p,\lambda;\C_r}\qquad \forall \ \mathbf{f}\in
L^{p,\lambda}(\C_r)
\end{align}
for any parabolic cylinder $\C_r\subset \ol{Q}.$
\end{lemma}
The bound \eqref{eqKf} has been proved by Fabes and Rivi\`ere
\cite[Theorem~1]{FR}
for general kernels of mixed homogeneity. In the case $b=1,$
\eqref{eqCaf} is obtained in \cite[Theorem~2.12]{BC1} while \cite{BC2} deals
with commutators of singular integrals with {\it constant} kernels
on homogeneous spaces. (Let us note that $\R^{n+1}$ endowed with the metric
$\rho(x,t)$ is a homogeneous space.) The passage from {\it constant} to  {\it
variable} kernels is standard, it makes use of Calder\'on-Zygmund's
approach (\cite{CZ1,CZ2}) of expansion
into spherical harmonics and leads to \eqref{eqCaf} for $b>1$ as well.

The estimates \eqref{eqKfMor}, \eqref{eqCafMor} have been obtained in authors'
paper \cite[Theorem~2.1, Corollary~2.7]{PS}
in the general case of kernels with mixed homogeneity.
Note that an adaptation of the proofs of \cite{PS} to the $L^p$-framework
gives \eqref{eqKf} and	\eqref{eqCaf}.

The estimates \eqref{3.10} and \eqref{3.11} follow from
\eqref{eqCaf} and \eqref{eqCafMor} on the base of $\A_{\alpha}\in VMO$ (see
\cite[Theorem~2.13]{CFL}, \cite[Corollary~2.8]{PS}, \cite[Theorem~3.7]{Sf}).

\begin{remark}\label{solution}\em
Employing density arguments, \eqref{eqKf} and \eqref{eqCaf},
it is easily seen that
the representation formula \eqref{10} still holds true (almost everywhere) for
 {\it compactly supported in $x$ functions\/} $\bv \in W^{2b,1}_p$ such that
$\bv(x,0)=\mathbf{0}.$
\end{remark}

\section{Remarks and Counterexamples}\label{sect5}
\setcounter{equation}{0}

\subsubsection{Parabolic systems with lower order terms.}
All the results presented here could be extended, modulo unessential
technicalities, to parabolic systems with lower order terms
\[
D_t \bu -
\sum_{|\alpha|=2b}\A_\alpha (x,t) D^\alpha \bu
+\sum_{|\beta|\leq 2b-1}\B_\beta (x,t) D^\beta \bu
=\mathbf{f}(x,t)
\]
with $\B_\beta(x,t)=\left\{b_\beta^{kj}(x,t) \right\}_{k,j=1}^m.$
Indeed, the coefficients $b_\beta^{kj}$  have to belong to suitable
Lebesgue
(Morrey) spaces chosen in such a way that $\B_\beta(x,t) D^\beta \bu$ and the
right-hand side $\mathbf{f}$ stay at the same space for any $\bu \in
W^{2b,1}_{p,\mathrm{loc}}(Q)$ ($\bu\in W^{2b,1}_{p,\lambda,\mathrm{loc}}(Q)$
respectively).

\subsubsection{Non-zero initial conditions.}
From a physical point of view, \eqref{system} governs a forward deterministic
process which evolves from its initial state. That is why
our  results  regard
strong solutions of the system \eqref{system} with {\it zero initial trace\/}
($\bu(x,0)=\mathbf{0}$). This is an important requirement ensuring us to deal
only with Gaussian-type (volume) potential
\begin{equation}\label{gauss}
\bu(x,t)=\int_{\R^{n+1}}\bs{\Gamma}_0(x-y,t-\tau)\LL \bu(y,\tau)dyd\tau
\end{equation}
in order to  derive the representation formula \eqref{10} for the highest order derivatives
$D^{2b}\bu.$
If $\bu(x,0)=\bs{\varphi}(x)$ then, roughly speaking, one should add a
{\it surface\/} potential
\[
\int_{\R^{n}} \bs{\Gamma}_0(x-y,t) \bs{\psi}(y) dy
\]
 to the right-hand side of \eqref{gauss},  where the density $\bs{\psi}$ depends of
$\bs{\varphi}(x)$ and $\mathbf{f}(x,t)$ (see \cite{Ed1,Ed2,F,LSU,Sl1} for details).
It is clear that the representation formula
\eqref{10} will change and this will reflect to the a'priori estimates proved in
Theorems~\ref{th1} and \ref{th3} by adding the respective norm of $\bs{\varphi}$
on $\Omega''$ to the right-hand sides of \eqref{9} and \eqref{9a}. (Note that the
initial data $\bs\varphi$ influence on the solution
$\bu(x,t)$ at {\it any\/} instant $t>0$ because  {\it linear\/} parabolic
systems support {\it infinite propagation speed of disturbances.\/})

To obtain the exact dependence on ${\bs\varphi}$ at \eqref{9} and \eqref{9a},
we follow an equivalent but more simple approach. In effect, the initial trace $\bs{\varphi}$
belongs to the Besov space $B^{2b-2b/p}_{p,\mathrm{loc}}(\Omega)$ (cf.
\cite[Theorem~5.1]{Sl1}) in
case of Theorem~\ref{th1}, while $\bs{\varphi}$ is a trace on $\{t=0\}$ of a
function lying in Sobolev--Morrey  class $W^{2b,1}_{p,\lambda}(Q)$ (cf. \cite{Cm0}
in the situation studied by Theorem~\ref{th3}). Therefore, without loss of
generality, we may suppose that $\bs{\varphi}(x)$ is extended (first in
$\R^n$ according
to Hestens--Whitney, and then in $\R^n\times(0,T)$ by solving Cauchy problem for
suitable system
with
{\it constant\/} coefficients) to a function $\bs{\Phi}(x,t)\in
W^{2b,1}_{p,\mathrm{loc}}(Q)$
($\bs{\Phi}\in W^{2b,1}_{p,\lambda,\mathrm{loc}}(Q)$) such that
$\bs{\Phi}(x,0)=\bs{\varphi}(x)$
for a.a. $x\in\Omega.$ Now, if $\bu$ solves \eqref{system} and
$\bu(x,0)=\bs{\varphi}(x)$
then $\tl\bu:=\bu-\bs{\Phi}$ is such that
\[
\LL\tl\bu=\mathbf{f}-\LL\bs{\Phi}\quad \text{a.e.}\ Q,\qquad
\tl\bu(x,0)=\mathbf{0}\quad \text{a.a.}\ x\in\Omega
\]
and therefore \eqref{9} (respectively, \eqref{9a}) holds for it. In other words,
$\bu$ will satisfy an a'priori bound
like \eqref{9} (or \eqref{9a}) with the norm
$\|\bs{\Phi}\|_{W^{2b,1}_{p}(Q'')}\sim
\|\bs{\varphi}\|_{B_p^{2b-2b/p}(\Omega'')}$ (respectively,
$\|\bs{\Phi}\|_{W^{2b,1}_{p,\lambda}(Q'')}$)
added to the right-hand side.

\subsubsection{Non-parabolic systems.} We are going to show that the parabolicity
condition \eqref{parabolicity} is {\it necessary\/} in order to be satisfied the a'priori estimate
\eqref{9}.

Consider the system with constant
coefficients
\begin{equation}\label{non}
\NP\bu:= (D_t-\ER(D_x))\bu=
D_t\bu - \left(
\begin{matrix}
\delta_1\Delta & 0 & \ldots & 0\\
0 & \delta_2\Delta &  \ldots & 0\\
\vdots & \vdots & \ddots & \vdots\\
0 & 0 & \ldots & \delta_m\Delta
\end{matrix}
\right)\bu=\mathbf{f}
\end{equation}
where $\delta_j\in\R$ and $\Delta$ is the Laplace operator.
 Let us emphasize that, in
contrast to the case of a single second-order equation ($m=b=1$), the operator $\NP$
is {\it not necessarily parabolic\/} even if the
matrix operator $\ER(D_x)$ is {\it strongly elliptic\/} in the sense that
\begin{equation}\label{r2}
\mathrm{det\,}\big(\ER(\xi)\big)\geq
K|\xi|^{2bm},\quad K=\mathrm{const}>0
\end{equation}
(cf. \cite{A,ADN,CFF,DN}). In fact, \eqref{r2}
reads $\prod_{j=1}^m \delta_j\geq K>0$ whereas $\NP$ is a parabolic operator
iff $\delta_j>0$ for each $j\in\{1,2,\ldots,m\}$ as it follows from \eqref{parabolicity}.

To show that, in general, the estimate \eqref{9} {\it cannot hold for
non-parabolic systems,\/} define $\Omega'=[-\pi,\pi]^n,$ $\Omega''=[-2\pi,2\pi]^n,$
$Q'=\Omega'\times(0,1),$ $Q''=\Omega''\times(0,1).$\\[4pt]
\indent
{\it Roots of \eqref{parabolic} with zero real part.\/}
Let $\delta_1=0.$ For each integer $N,$
$\bu_N(x,t):=\big(t\sin(N x_1),0,\ldots,0\big)^{\mathrm{T}}$ solves
\eqref{non} with $\mathbf{f}(x,t)=\big(\sin(N x_1),0,
\ldots,0\big)^{\mathrm{T}}$ and $\bu_N(x,0)=\mathbf{0}.$
Straightforward calculations give $\|\bu_N\|_{W^{2,1}_2(Q')}\sim N^2(1+o(1))$
as $N\to +\infty$ whereas $\|\bu_N\|_{L^{2}(Q'')}+\|\mathbf{f}\|_{L^{2}(Q'')}
=\mathrm{const}.$ Hence \eqref{9} {\it cannot be satisfied\/} by
$\bu_N$ for large $N.$\\[6pt]
\indent
{\it Roots of \eqref{parabolic} with positive real part.\/}
Consider \eqref{non} with $\delta_1=-1.$ The function
\[
\bu_N(x,t):=\big(({\mathrm{e}}^{N^2 t}-1)\sin(N x_1),0,\ldots,0\big)^{\mathrm{T}}\qquad
\forall N\in\N
\]
is a solution of \eqref{non} with $\mathbf{f}(x,t)=\big(N^2\sin(N x_1),0,
\ldots,0\big)^{\mathrm{T}}$ and $\bu_N(x,0)=\mathbf{0}.$
 However,
\[
\|\bu_N\|_{W^{2,1}_2(Q')}\sim N {\mathrm{e}}^{N^2}(1+o(1)),\quad
\|\bu_N\|_{L^{2}(Q'')}+\|\mathbf{f}\|_{L^{2}(Q'')} \sim N^{-1} {\mathrm{e}}^{N^2}(1+o(1))
\]
as $N\to +\infty,$ and therefore \eqref{9}
{\it fails\/} once again for $\bu_N$ when $N$ is large enough.
The first component $(u_{N})_1(x,t)$ of $\bu_N$ solves the  {\it
backward\/} heat equation
\begin{equation}\label{backward}
\big(D_t+\Delta\big) (u_{N})_1=(f_{N})_1\quad \text{in}\ Q''
\end{equation}
for which the Cauchy problem with initial data given at $t=0$ is
{\it not well-posed\/} in sense of Hadamard (the process $(u_{N})_1(x,t)$
cannot be recovered from its {\it final\/} state $(u_{N})_1(x,0)$!).
Actually, reversing the time
$\tau=1-t$ and setting $\mathbf{U}_N(x,\tau)=\bu_N(x,1-\tau)$ for
$\tau\in[0,1],$
it turns out that $(\mathbf{U}_N)_1(x,\tau)$ solves the {\it forward\/}
heat equation
\[
\big(D_\tau-\Delta\big) (U_{N})_1=-(f_{N})_1\quad \text{in}\ Q''
\]
for which $(U_N)_1(x,\tau)|_{\tau=0}$ and $(U_N)_1(x,\tau)|_{\tau=1}$ are
the {\it initial\/} and {\it final\/} state, respectively. The failure of
Theorem~\ref{th1} for $\bu_N(x,t)$ is due, therefore, to
initial data given at wrong instant --- these are to be taken at $t=1.$ In fact,
the function
$\bv_N(x,t):=\big(({\mathrm{e}}^{N^2 t}-{\mathrm{e}}^{N^2})\sin(N x_1),0,\ldots,0\big)^{\mathrm{T}}$ is a
solution of
\eqref{non} with $\delta_1=-1$ and
$\mathbf{f}(x,t)=\big(N^2{\mathrm{e}}^{N^2}\sin(N x_1),0,\ldots,0\big)^{\mathrm{T}},$
and it satisfies
\eqref{9}. Note that $\bv$ has
zero trace on the plane $\{t=1\}.$

Turning back to \eqref{non}, suppose $\delta_j<0$ for {\it each\/} $j\in\{1,2,\ldots,m\}.$
Then
$\NP$ is not a parabolic operator, but a simple reversal of time as above
reduces
\eqref{non} to an equivalent parabolic system for which the initial data are
to given
at $t=T.$ This procedure fails categorically when \eqref{non} contains
{\it both\/} positive and negative $\delta_j$ and then also Theorem~\ref{th1}
collapses
for $\NP.$ In fact, for a fixed $j\in\{1,2,\ldots,m\}$ the initial trace of the
component
$u_j(x,t)$ must be imposed at $t=0$ if $\delta_j>0$ or at $t=T$ if
$\delta_j<0$ and
this contradicts to the evolutionary nature of the system.

\subsubsection{Counterexamples} Estimates similar to our local bounds \eqref{6} have been
announced recently by M.A. Ragusa in \cite{Ra}.
The author deals with systems of the type
\begin{equation}\label{r1}
D_t \bu - \ER(x,t,D_x)\bu +\BR(x,t,D_x)\bu=\mathbf{f}(x,t)
\end{equation}
where $\ER(x,t,D)=\big\{\sum_{|\alpha|=2b}a^{(\alpha)}_{jk}(x,t)
D^\alpha\big\}_{j,k=1}^m$ and $\BR(x,t,D)$ is a linear matrix-differential operator
of
order less than $2b.$ The principal coefficients $a^{(\alpha)}_{jk}$ belong to
$VMO\cap L^\infty$ while these of $\BR$ are in $L^q$ with some $q,$ and it is
supposed that $\ER(x,t,D)$ is a {\it strongly elliptic\/} operator.

Setting $Q_\sigma= \{(x,t)\in\R^{n+1}\colon\ |x-x_0|<\sigma,\ |t-t_0|<\sigma^2\},$
the author asserts
\begin{equation}\label{r3}
\|D^{2b}\bu\|_{L^p(Q_{\sigma/2})}+\|D_{t}\bu\|_{L^p(Q_{\sigma/2})}\leq
C\big(\|\mathbf{f}\|_{L^p(Q_{\sigma})}+\|\bu\|_{L^p(Q_{\sigma})}\big)
\end{equation}
(compare with \eqref{6} and note the factor $Cr^{-2b}$ of $\|\bu\|_{p,r}$
there) for any	solution  of \eqref{r1} and any $Q_\sigma\Subset Q,$ and
\begin{equation}\label{r4}
\|D^{2b}\bu\|_{C^{0,\gamma}(Q_{\sigma/2})}\leq
C\big(\|\mathbf{f}\|_{L^p(Q_{\sigma})}+\|\bu\|_{L^p(Q_{\sigma})}\big)
\quad\text{with}\ \gamma=1-{(n+1)}/{p}.
\end{equation}

As \eqref{non} shows however, \eqref{r1} is {\it not,\/} in general, a
{\it parabolic\/} system. Thus, \eqref{r2} means that the characteristic equation
\eqref{parabolic} corresponding to \eqref{r1}
could have roots with {\it positive\/} real part and then the contour
$\mathfrak{C}(\xi)$ in
\eqref{fund-equa} will necessarily pass through the right complex half-plane
causing this way loss of analyticity of the integrand in \eqref{fund-equa}
(cf. \cite[Appendix~I]{Sl1}).
In effect, instead of
\eqref{det} the author considers the equation
$u_t-\text{det\,}(\sum_{|\alpha|= 2b}a^{(\alpha)}_{jk}D^\alpha)u=0$
which is not parabolic in general,
and therefore its fundamental solution cannot be analytic as affirmed.

Regarding the estimate \eqref{r3}, it is rough even if
\eqref{r1} was a parabolic system.
In fact, \eqref{r3} provides for a control of $D^{2b}\bu(x,t)$ and
$D_t\bu(x,t)$ for
$t\in(t_0-\sigma^2/4,t_0+\sigma^2/4)$ in terms of the future states $\bu(x,t)$
with
$t\in(t_0+\sigma^2/4,t_0+\sigma^2)$ (compare the $b$-independent $Q_\sigma$'s
with the cylinders \eqref{pc} which depend on $b$ and  have the same
upper base). Moreover, $Q_\sigma\Subset Q$ and this presumes
the upper base $\{(x,T)\colon\ x\in\Omega\}$ of $Q$ to be a part of the
boundary which
is not the case when dealing with parabolic operators. It is difficult therefore
to imagine
how a global-in-time estimate like \eqref{9} could be derived from \eqref{r3}.

As shown above, the mixed homogeneity of degree $-n-2b$ of the
derivatives
$D^{2b}_y \bs{\Gamma}(x,t;y,\tau)$ (cf. $(\P_2)$) is a crucial property
ensuring
these are kernels of Calder\'on--Zygmund type, and therefore
validity of Lemma~\ref{pr1}. Instead, the claim in \cite{Ra} is that the
$2bm$-order
derivatives of the fundamental solution  of \eqref{r1} (that is,
$D^{2b}_y\bs{\Gamma}(x,t;y,\tau)$ when \eqref{r1} is a parabolic system)
are homogeneous of degree $-n$ which is
false even for the single heat equation.

As consequence of all these gaps, also the H\"older continuity \eqref{r4} is
wrong as we are
going to show on the level of a very classical example (see the Introductory
Chapter in \cite{LSU}).
Let $p>n+1,$ $Q=\big\{(x,t)\in \R^n\times\R^+\colon\ |x|\leq 1,\
t\in(0,2)\big\},$
and consider the function $v(x,t)=|x|^{2\mu}\frac{|x|^{2}}{|t-1|}
\exp\left\{-\frac{|x|^{2}}{|t-1|}\right\}$
with $\mu\in \left(1-\frac{n+2}{2p},1-\frac{n+1}{2p}\right)\subset(0,1).$
Straightforward analysis yields $D_t v-\Delta v= f(x,t)$ a.e. in $Q,$
$f\in L^p(Q),$ and $|D_xv(x,t)|=|x|^{2\mu-1}(1+O(1))$ as $(x,t)\to(0,1).$

Let $p\geq n+2.$ Then $1-\frac{n+2}{2p}\geq \frac{1}{2}$ and $2\mu-1\in
(0,1).$ If
$p\in(n+1,n+2)$ then $1-\frac{n+2}{2p}<\frac{1}{2}$ and restricting $\mu$ to
the interval
$\left(\frac{1}{2},1-\frac{n+1}{2p}\right)$ we get $2\mu-1\in(0,1)$ once
again. In both cases,
the gradient $D_xv(x,t)$ is H\"older continuous with exponent $2\mu-1$ near
the point $(0,1),$
but $2\mu-1<1-\frac{n+1}{p}$ and therefore $D_xv$ {\it cannot\/} be H\"older
continuous with
the exponent $\gamma$ as \eqref{r4} asserts. Moreover, taking $\mu\in
\left(1-\frac{n+2}{2p},\frac{1}{2}\right)$ when $p\in (n+1,n+2),$ $2\mu-1$
becomes negative
and $|D_xv(x,t)|$ even blows-up at $(0,1)$ contrary to the statement
\eqref{r4}.

\vspace{0.2pc}

In view of the above comments, the announcements of \cite{Ra} can hardly be
considered proven with the necessary
rigour.


\begin{thebibliography}{CFL23}
\bibitem{A}
{ S. Agmon},
{\em Lectures on Elliptic Boundary Value Problems} (Van Nostrand, Princeton,
N.J., 1965).
\bibitem{ADN}
{ S. Agmon, A. Douglis \and L. Nirenberg},
`Estimates near the boundary for solutions of elliptic partial
differential equations satisfying general boundary conditions', I,
{\em Commun.\ Pure\ Appl.\ Math.\ }{12} (1959) 623--727; II, ibid. {17} (1964)
35--92.
\bibitem{BC1}
{ M. Bramanti \and M.C. Cerutti},
`$W^{2,1}_p$-solvability
for the Cauchy--Dirichlet problem for parabolic equations with $VMO$
coefficients',
{\em Comm.\ Part.\ Diff.\ Equations\ }{18} (1993)
1735--1763.
\bibitem{BC2}
{ M. Bramanti \and M.C.  Cerutti},
`Commutators of singular integrals on homogeneous spaces',
{\em Boll.\ Un.\ Mat.\ Ital.\ B\ (VII)\ }{10} (1996)
843--883.
\bibitem{CZ1}
{ A.P. Calder\'on \and A. Zygmund},
`On the existence of certain singular integrals',
{\em Acta\  Math.\ }{88}  (1952)  85--139.
\bibitem{CZ2}
{ A.P. Calder\'on \and A. Zygmund},
`Singular integral operators and differential equations',
{\em Amer.\ J.\ Math.\ }{79}  (1957)  901--921.
\bibitem{Cm0}
{ S. Campanato},
`Caratterizzazione delle tracce di funzioni appartenenti ad
una classe di Morrey insieme con le loro derivate prime',
{\em Ann.\ Scuola\ Norm.\ Sup.\ Pisa }{15} (1961) 263--281.
\bibitem{Cm}
{ S. Campanato},
`Propriet\'a di H\"olderianit\`a di alcune classi
di funzioni',
{\em Ann.\ Scuola\ Norm.\ Sup.\ Pisa }{17} (1963) 175-188.
\bibitem{Cm2}
{ S. Campanato},
{\em Sistemi ellittici in forma divergenza. Regolarit\`a
all'interno,}
(Pubblicazioni della Classe di Scienze: Quaderni,
Scuola Norm. Sup., Pisa, 1980).
\bibitem{CFF}
{ F. Chiarenza, M.~Franciosi \and  M.~Frasca},
`$L^p$ estimates for linear
elliptic systems with discontinuous coefficients',
{\em Rend.\ Accad.\ Naz.\ Lincei,\ Mat.\ Appl.\ }{5} (1994) 27--32.
\bibitem{CFL}
{ F. Chiarenza,  M. Frasca \and P.~Longo},
`Interior $W^{2,\,p}$ estimates for non divergence elliptic
equations with discontinuous coefficients',
{\em Ric.\ Mat.\ }{60} (1991) 149--168.\bibitem{Da}
{ G. Da Prato},
`Spazi ${\mathcal{L}}^{p,\theta}(\Omega,\delta)$ e loro propriet\`a',
{\em Ann.\ Mat.\ Pura Appl.\ (IV)\ }{69} (1965) 383--392.
\bibitem{DN}
{ A. Douglis \and L. Nirenberg},
`Interior estimates for elliptic systems of partial differential equations',
{\em Commun.\ Pure\ Appl.\ Math.\ }{8} (1955) 503--538.
\bibitem{Ed1}
{ S.D. Ejdel'man},
{\em Parabolic\ Systems,} (North-Holland, Amsterdam, London, 1969).
\bibitem{Ed2}
{ S.D. Ejdel'man},
`Parabolic Equations', In: Yu. V. Egorov and M. A. Shubin (Eds.) {\em Partial\
Differential\ Equations\ VI:\ Elliptic\ and\ Parabolic\ Operators,}
Encycl. Math. Sci., \textbf{63}, pp. 203--316, (Springer-Verlag, Berlin,
1994).
\bibitem{FR}
{ E.B.~Fabes   \and  N.~Rivi\`ere},
`Singular integrals with mixed homogeneity',
{\em Studia\ Math.\ }{27} (1966) 19--38.
\bibitem{F-ind}
{ A. Friedman},
`Interior estimates for parabolic systems of partial differential equations',
{\em J.\ Math.\ and\ Mech.\ }{7} (1958) 393--418.
\bibitem{F}
{ A. Friedman},
{\em Partial\ Differential\ Equations\ of\ Parabolic\ Type,}
(Prentice-Hall, Englewood Cliffs, N.J., 1964).
\bibitem{GS}
{ I.M. Gel'fand  \and  G.E. Shilov},
{\em Generalized\ Functions. Vol.~2: Spaces of
Fundamental and Generalized Functions,}
(Academic Press, New York, London, 1968).
\bibitem{I}
{ V.P. Il'in},
`Properties of certain classes of differentiable functions of
several variables defined in an $n$-dimensional domain',
{\em Trudy\ Mat.\ Inst.\ Steklov\ }{66} (1962) 227--363.
\bibitem{ISl}
{ V.P. Il'in \and  V.A. Solonnikov},
`Some properties of differentiable functions of several variables',
{\em Trudy\ Mat.\ Inst.\ Steklov\ }{66} (1962) 205--226.
\bibitem{JN}
{ F.~John  \and  L.~Nirenberg},
`On functions of bounded mean oscillation',
{\em Commun.\ Pure\ Appl.\ Math.\ }{14} (1961) 415--426.
\bibitem{K}
{ M. Kronz},
`Some function spaces on spaces of homogeneous type',
{\em Manuscr.\ Math.\ }{106} (2001) 219-248.
\bibitem{LSU}
{ O.A.~Ladyzhenskaya, V.A.~Solonnikov \and  N.N.~Ural'tseva},
{\em Linear\ and\ Quasilinear\ Equations\ of\ Parabolic\ Type}, Transl. Math.
    Monographs, \textbf{23},  (Amer. Math. Soc., Providence, R.I., 1968).
\bibitem{MPS}
{ A. Maugeri, D.K. Palagachev \and L.G. Softova},
{\em Elliptic\ and\ Parabolic\ Equations\ with\ Discontinuous\ Coefficients,}
(Wiley-VCH, Berlin, 2000).
\bibitem{M}
{ C.B. Morrey Jr.},
{\em Multiple\ Integrals\ in\ the\ Calculus\ of\
Variations,} (Springer-Verlag, Berlin, 1966).
\bibitem{PS}
{ D. Palagachev \and  L. Softova},
`Singular integral operators, Morrey spaces and fine
regularity of solutions to PDE's',
{\em Potential\ Anal.\ }{20} (2004) 237-263.
\bibitem{Ra}
{ M.A. Ragusa},
`Parabolic systems with non-continuous coefficients',
In: {\em Proc.\ Conf.\ on\ Dynamical\ Systems\ and\ Differential\ Equations,
(Wilmington, NC, 2002)},
{\em Discrete\ Contin.\ Dyn.\ Syst.\ } suppl. (2003) 727--733.
\bibitem{S}
{ D. Sarason},
`Functions of vanishing mean oscillation',
{\em  Trans.\ Amer.\ Math.\ Soc.\ }{207} (1975) 391--405.
\bibitem{Sf}
{ L. G. Softova},
`Parabolic equations with $VMO$ coefficients in Morrey spaces',
{\em  Electr.\ J.\ Diff.\ Equations\ }{2001} (2001), Paper No. 51.
\bibitem{Sl1}
{ V. A. Solonnikov},
`On the boundary value problems for linear parabolic systems of differential
equations of general form',
{\em Proc.\ Steklov\ Inst.\ Math.\ }{83} (1965);
English translation: O. A. Ladyzhenskaya~(Ed.)
{\em Boundary\ Value\ Problems\ of\ Mathematical\ Physics\ III,}
(Amer. Math. Soc., Providence, R.I., 1967).
\bibitem{Sl2}
{ V. A. Solonnikov},
`Estimates in $L_p$ of solutions of elliptic and parabolic systems',
{\em Proc.\ Steklov\ Inst.\ Math.\ }{102} (1967) 137--160;
English translation: O. A. Ladyzhenskaya~(Ed.)
{\em Boundary\ Value\ Problems\ of\ Mathematical\ Physics\ V,}
(Amer. Math. Soc., Providence, R.I., 1970).
\end{thebibliography}
\end{document}